\documentclass{article}

\usepackage{arxiv}

\usepackage[utf8]{inputenc} % allow utf-8 input
\usepackage[T1]{fontenc}    % use 8-bit T1 fonts
\usepackage{hyperref}       % hyperlinks
\usepackage{url}            % simple URL typesetting
\usepackage{booktabs}       % professional-quality tables
\usepackage{amsfonts}       % blackboard math symbols
\usepackage{nicefrac}       % compact symbols for 1/2, etc.
\usepackage{microtype}      % microtypography
\usepackage{lipsum}		% Can be removed after putting your text content
\usepackage{graphicx}
\usepackage{natbib}
\usepackage{doi}
\usepackage[linesnumbered,ruled]{algorithm2e}

\SetKw{KwBy}{by} 
\usepackage{amsmath, amssymb}
\usepackage{enumitem}
\usepackage{booktabs}
\usepackage{geometry}
\title{Assignment–Routing Optimization : Efficient Heuristic Solver with Shaking Algorithm}

\author{
  {\includegraphics[scale=0.06]{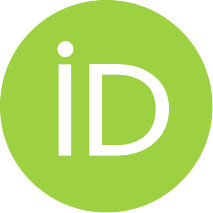}\hspace{1mm} Yuan Qilong $^{*a}$,  Michal Pavelka$^{b}$. } \\$^{a}$ Singapore Institute of Technology. $^{b}$ Mathematical Institute, Faculty of Mathematics and Physics, Charles University
\\  \texttt{$^*$ Corresponding Author Email: qilong.yuan@singaporetech.edu.sg} \\
}

\hypersetup{
  pdftitle={Assignment–Routing Optimization : Efficient Heuristic Solver with Shaking Algorithm},
  pdfsubject={Optimization, Combinatorics},
  pdfauthor={Yuan Qilong, Michal Pavelka},
  pdfkeywords={Joint Routing–Assignment Optimization, Heuristic Solver, Shaking Algorithm},
}

\begin{document}
\maketitle

\begin{abstract}
This paper works on heuristic solver for joint assignment and routing optimization problem. Study on previous works shows that MIP based exact solvers can only provide efficient solutions for small to moderate size problems, due to exponentially growing computational complexity. This paper proposes to start with high quality initial guess through Hungarian algorithm based assignment and  heuristic cycle merging algorithm. Subsequently, the solution is improved based on a proposed shaking algorithm to improve the assignment and routing sequence. In addition, the shaking approach also enables the Simulated Annealing algorithm to further improve the solution, which is very difficult if it is purely based on random sampling updates of item and placeholder sequences. Extensive experimental validation comparing with ground truth from the previously shared database shows that the introduced solver is much more efficient than the Gurobi solver especially for large size problems, with a 1000 node pair problem being solved within 1 min in Python implementation. The solution accuracy is within a percent in general as compared with ground truth in database. Although there are spaces for the proposed solver to be further improved with better accuracy, it works for practical applications with acceptable path quality and sufficient solver efficiency. Such shaking algorihtms based solvers can also be applied to more general joint assignment and routing optimization problem with multiple type of items and corresponding placeholders. GitHub repository: \url{https://github.com/QL-YUAN/Joint-Assignment-Routing-Optimization-Heuristic.git}
\end{abstract}

% keywords can be removed
\section{Introduction}
\label{sec:intro}

Many real-world applications in robotics, logistics, and industrial automation require transporting a set of discrete items from their initial positions to designated destinations. These tasks naturally intertwine three interdependent decision layers: 
(i) \emph{assignment} — deciding which object should go to which target location; 
(ii) \emph{sequencing} — determining the order of pick-and-place actions; and 
(iii) \emph{routing or motion planning} — computing feasible paths that minimize travel cost under geometric and dynamic constraints. This paper focuses on the first two layers, while assuming that underlying motion planning can be handled separately or abstracted into travel cost estimates.

Such problems arise in a variety of domains:

\begin{itemize}
  \item \textbf{Vehicle Routing Problem with Pickup and Delivery(VRPPD).} 
  This problem involves mobile robots or autonomous vehicles tasked with transporting items or passengers between specific locations.  General VRPPD problems are of different variations with high complexity \cite{parragh2008survey}. The vehicle routing and traveling salesman problem with pickup and delivery (TSPPD, VRP-PD) involve matching pickup–delivery pairs under precedence and capacity constraints \cite{HERNANDEZPEREZ2004126}. The stacker-crane abstraction \cite{hernandez2004,hernandez2007,treleaven2012} captures a related setting where a single vehicle transports all items between paired locations, forming a constrained Hamiltonian tour. These problems are NP-hard and have motivated both exact methods (e.g., branch-and-cut) and scalable heuristics. In case of one vehicle having unit capacity to pick and deliver one customer at a time with a routing mode, this becomes one variation of Travel Salesman Problem. However, if the picking and delivery pair becomes flexible, the problem becomes more complex.   
  \item \textbf{Tabletop rearrangement.} Robotic manipulators tasked with relocating scattered objects (e.g., puzzle pieces, tools) into designated slots or target poses, where assignment and action sequencing are critical to efficiency. Hierarchical task-and-motion planning (TAMP) strategies interleave discrete decision-making (e.g., which object to move) with motion planning \cite{simeon2004, king2016rearrangement}, for pick and placing items with optimal robot actions \cite{han2019high}.
\end{itemize}

Based on the application requests, in this work, the authors specifically study the optimization problem to find the near-optimal path that starts from a depot and iteratively visits all possible items-to-placeholder pairs and finally goes back to the depot, forming a Hamiltonian cycle. 

\subsection{Existing Solvers}
Despite progress in assignment and routing subproblems, solving their integrated form—particularly under precedence, capacity, and connectivity constraints—remains a fundamental challenge.
While state-of-the-art TSP solvers (e.g., Concorde \cite{applegate2006traveling}) are highly optimized for standard symmetric TSPs, they struggle when additional constraints—such as precedence, unit capacity, or unknown assignments—are introduced. Attempts to encode the problem into a standard TSP using gadgetization often result in graph blow-up and loss of exploitable structure, and constraints for item and placeholder pairs is likely to break based on experimental testing. Similarly, practical routing libraries like Google OR-Tools \cite{ortools2024} provide heuristics for pickup and delivery, but typically assume fixed pairings or require large numbers of additional assignment variables, limiting scalability\cite{hernandez2004,hernandez2007,treleaven2012}.
\subsubsection{MIP solvers and open challenges}

In previous work, the authors introduces a \textbf{Mixed-Integer Programming (MIP)} formulation that jointly encodes object-to-goal assignments, pickup-before-delivery precedence, and Hamiltonian tour constraints\cite{yuan2025}. This model is solved using the Gurobi solver, enhanced with dynamic cutting-plane callbacks to eliminate subtours efficiently. The \textbf{benchmark datasets} and corresponding \textbf{ground-truth optimal solutions} are also released, which serve both as a reproducible evaluation platform and a reference for heuristic algorithm development.

Experimental results show that the exact MIP solver performs effectively for small to moderate-sized instances (up to 100 pickup--delivery pairs, solving time is within 1 seconds), consistently finding optimal solutions in reasonable time. However, as the problem scales beyond 300 pairs, the computational cost increases rapidly, and Gurobi can fail to return a feasible solution within acceptable time limits when $n$ reaches 1000. This reinforces known theoretical complexity results and highlights the practical limitations of exact solvers for large-scale joint assignment--routing problems.

As a result, there is clear motivation for developing efficient heuristics or approximation algorithms that can deliver near-optimal solutions quickly, particularly when responsiveness is essential. 

\subsection{Highlighs of This Work}

The highlights of this work include:

\begin{enumerate}
  \item A novel, efficient heuristic that combines Hungarian \cite{kuhn1955hungarian} assignment, with a proposed shaking algorithm to  improve the assignment and the precedence efficiently, followed by  simulated annealing–based sequencing (and shaking), achieving near-optimal results with significantly reduced runtime.
  \item A comparative study was conducted across exact (Gurobi), heuristic (the introduced method), and library-based (OR-Tools) approaches using synthetic benchmarks with up to 300 pickup–delivery pairs. The results demonstrate both the accuracy and scalability of the introduced heuristic. The method achieves a balance of speed and precision, remaining within 1\% of optimal for moderate-scale problems, scaling efficiently to 300 node pairs in under 5 seconds, and solving 1000 node pairs in approximately one minute. While there is room for improvement in solution quality to achieve even better paths, this heuristic solver provides practical and efficient solutions suitable for real-world applications.
  \item GitHub repository : \url{https://github.com/QL-YUAN/Joint-Assignment-Routing-Optimization}
\end{enumerate}

\section{Problem Formulation}
\subsection{Scenario Description}

Consider a set of $N = n$ items $\{I_1, I_2, \dots, I_n\}$ located arbitrarily in a 2D workspace, and a corresponding set of $N$ placeholder locations $\{S_1, S_2, \dots, S_n\}$.

\textbf{Assignment:} Each item must be assigned to a unique placeholder (one-to-one) and delivered there by a single agent.

The agent begins at an initial position $p_\text{start} \in \mathbb{R}^2$, picks up one item at a time, and immediately delivers it to a uniquely assigned placeholder, successively visiting each item $I_1, \dots, I_n$ in a specific order. After completing all deliveries, the agent returns to $p_\text{end}$. If $p_\text{start} = p_\text{end}$, the path forms a closed loop.  

To present the formulation compactly, $p_\text{end}$ is assigned to $I_0$, and $p_\text{start}$ is assigned to $S_0$, with a fixed assignment from $I_0$ to $S_0$, ensuring that a closed loop can be formed.  

\subsection{Notation}

\begin{itemize}[leftmargin=1.5em]
  \item Let $I = \{0, 1, \dots, n\}$ be the set of item indices, where $I_0$ corresponds to the ending position.
  \item Let $S = \{0, 1, \dots, n\}$ be the set of placeholder indices, where $S_0$ corresponds to the start position.
  \item For each $i \in I$, let $p_i \in \mathbb{R}^2$ be the position of item $I_i$, with $p_0 = p_\text{end}$.
  \item For each $j \in S$, let $s_j \in \mathbb{R}^2$ be the position of placeholder $S_j$, with $s_0 = p_\text{start}$.
  \item Let $d(u,v)$ denote the Euclidean distance between points $u, v \in \mathbb{R}^2$.
  \item Let $a : I \to S$ be a bijective assignment of items to placeholders, with $a(0) = 0$.
  \item Let $\pi = (\pi_0, \pi_1, \dots, \pi_n)$ be a permutation of $I$ representing the pickup order, with $\pi_0 = 0$ fixed.
\end{itemize}

\subsection{Optimization Objective}

The agent follows the sequence:
\[
p_{\pi_0} \to s_{a(\pi_0)} \to p_{\pi_1} \to s_{a(\pi_1)} \to \dots \to p_{\pi_n} \to s_{a(\pi_n)} \to p_{\pi_0}
\]
That is, the agent starts from \( p_{\pi_0} = p_0 \)
, picks up item \( \pi_1 \) at \( p_{\pi_1} \). Then, delivers it to \( s_{a(\pi_1)} \), and proceeds to the next pickup location \( p_{\pi_{2}} \)..., till picking item at \( p_{\pi_{2}}\) and placing item at \(s_{a(\pi_n)}\). Finally, returns to \( p_{\pi_0} = p_0 \).

Then, the total travel cost is:

\[
T(\pi, a) = \sum_{k=0}^{n}[ d\big(p_{\pi_k},\, s_{a(\pi_k)}\big) + d\big(s_{a(\pi_k)},\, p_{\pi_{(k+1) \bmod (n+1)}}\big) ]
\]

Where:
\[
\pi_{(k+1) \bmod (N+1)} =
\begin{cases} 
\pi_{k+1}, & 0 \le k < N,\\[2mm]
\pi_0, & k = N
\end{cases}
\]

The goal is to find a bijective assignment \( a \) and a pickup sequence \( \pi \) (with \( \pi_0 = 0 \)) that minimize the total travel cost:
\[
(a^*, \pi^*) = \arg\min_{a,\pi} \; T(\pi, a)
=
\arg\min \big\{ T(\pi, a) \;\big|\; 
a : I \to S \text{ bijection},\ a(0)=0,\ \pi \in \text{Perm}(I),\ \pi_0=0 \big\}
\]

\section{Problem Solving Methods}
\begin{figure}[h]
  \centering  \includegraphics[width=0.6\linewidth]{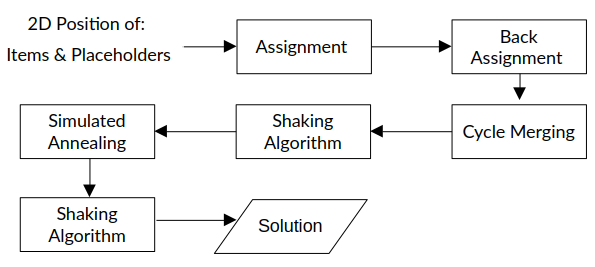}  
  \caption{Illustration of Workflow of The Proposed Hungarian algorithm}
  \label{fig:flow1}
  \vspace{0mm}
\end{figure}
In this work, as shown in Figure \ref{fig:flow1},  the joint assignment-routing problem is solved based on following steps:
\begin{itemize}
    \item Assign placeholders to items based on Hungarian Algorithm. 
    \item Back Assign items to different placeholders based on Hungarian Algorithm. 
    \item Create high quality initial solution based on cycle detection and cycle merging algorithms
    \item Efficiently improve the solution based on a proposed shaking algorithm. 
    \item Perform shaking algorithm enhanced Simulated Annealing Algorithm to improve the solution. 
    \item Final shaking algorithm can be added for further improvement.
\end{itemize}
\subsection{Item-Placeholder Assignment}
\begin{figure}[h]
  \centering  \includegraphics[width=0.6\linewidth]{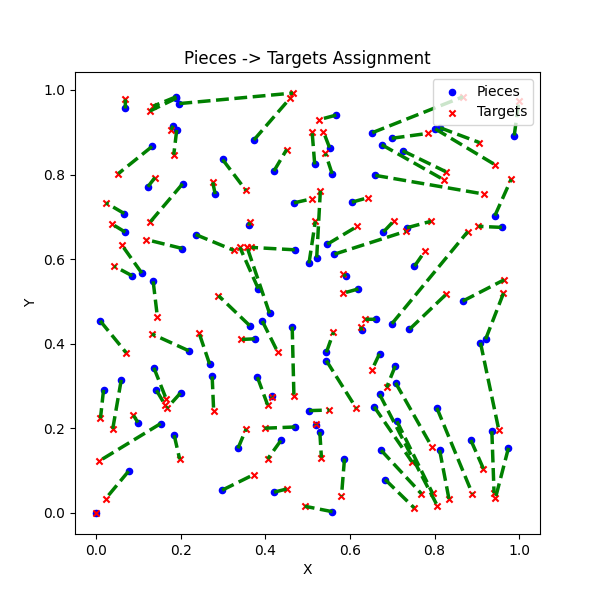}  
  \caption{Assignment of items towards placeholders with Hungarian algorithm}
  \label{fig:f1}
  \vspace{0mm}
\end{figure}
We first construct a cost matrix \( C^{(f)} \in \mathbb{R}^{(n+1) \times (n+1)} \), where each entry represents the Euclidean distance from item \( I_i \) to placeholder \( S_j \):
\[
C^{(f)}_{i,j} := d(p_i, s_j), \quad \text{for all } i,j \in \{0, \dots, n\}.
\]

We enforce the fixed assignment constraint:
\[
C^{(f)}_{0,j} := 
\begin{cases}
0 & \text{if } j = 0, \\
\infty & \text{otherwise}.
\end{cases}
\]

Then, we apply the \textbf{Hungarian algorithm} to find a bijective assignment:
\[
a^* : I \to S,
\]
minimizing the total assignment cost:
\[
a^* =
\arg\min_{a}
\left\{
\sum_{i=0}^{n} d(p_i, s_{a(i)})
\;\middle|\;
\begin{array}{l}
a : I \to S \text{ is a bijection},\\[2mm]
a(0) = 0
\end{array}
\right\}.
\]

\begin{figure}[h]
  \centering  \includegraphics[width=0.6\linewidth]{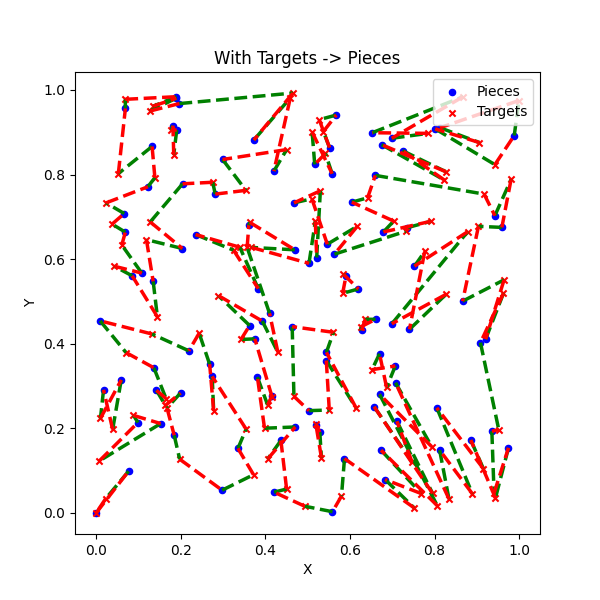}  
  \caption{Back assignment of items towards placeholders with Hungarian algorithm}
  \label{fig:f2}
\end{figure}

This yields a forward mapping \( a^* \) that assigns each item to a unique placeholder, satisfying all problem constraints including the fixed pair \( a^*(0) = 0 \). Figure \ref{fig:f2} provide sample of result for 100 pairs.

\subsection{Backward Connection: Placeholder-Item}
\vspace{-10mm}
After obtaining the forward assignment, we construct a second cost matrix \( C^{(b)} \in \mathbb{R}^{(n+1) \times (n+1)} \), where each entry \( C^{(b)}_{j,i} \) represents the cost of traveling from placeholder \( S_j \) to item \( I_i \):
\[
C^{(b)}_{j,i} := d(s_j, p_i), \quad \text{for all } j,i \in \{0, \dots, n\}.
\]

To prevent placeholder from linking back to the item it was originally assigned to, we define a masking operation:
\[
C^{(b)}_{a^*(i),\,i} := \infty, \quad \forall i \in I.
\]

This ensures that placeholder \( S_{a^*(i)} \) cannot link back to item \( I_i \) again in the loop construction.

Applying the Hungarian algorithm to \( C^{(b)} \), we compute a \textbf{backward assignment}:

\[
b^* = 
\arg\min_{b : S \to I}
\left\{
\sum_{j=0}^{n} d(s_j, p_{b(j)})
\;\middle|\;
\begin{array}{l}
b \text{ is a bijection}, \\[2mm]
b(a^*(i)) \neq i,\ \forall i \in I
\end{array}
\right\}.
\]

This two-step solution yields:
\begin{itemize}
    \item A forward delivery assignment \( a^* \), used for routing pickups to corresponding drop-offs,
    \item A back-connection mapping \( b^* \), useful for planning a cyclic return route, e.g., for routing heuristics.
\end{itemize}
Both steps rely on the Hungarian algorithm, which solves the assignment problem in \( \mathcal{O}(n^3) \) time.
The masking step in the backward phase ensures no trivial cycles are introduced. However, \textbf{if without marking to prevent placeholder from linking back to the item it was originally assigned to, it will result in a useful weak lower boundary of the path length.} Meanwhile, it is worth to mention that if a different $a \neq a^*$ function is chosen, the resulting b could leads to a different approximation solution, which worth to be further investigated for efficient evaluation algorithms. In this work, the initial guess starts from Hungarian assignment.

It is important to note that the combined forward and backward assignments do not guarantee a single-cycle solution. Specifically, the backward assignment constructed via the Hungarian algorithm may result in multiple disjoint cycles (i.e., a multi-cycle solution). This occurs because the backward connections are optimized independently of the forward paths, and no global constraint is enforced to merge all cycles into a single tour.

In applications where a single-loop traversal is required (e.g., for robotic execution or vehicle routing), further post-processing is necessary to ensure that the resulting plan is feasible for single loop-based execution. In the MIP solver, the cutting-plane based methods is computationally heavy mainly because of such cycle elimination process. Therefore, in this work, a cycle merging approximation method is applied. 

\subsection{Cycle Detection and Merging Procedure}

The two-phase Hungarian assignment procedure described above typically yields multiple disjoint cycles, rather than a single tour, as shown in Figure \ref{fig:f2}. To construct a feasible solution respecting the single-cycle constraint, we propose a cycle merging heuristic consisting of the following steps:

\paragraph{Step 1: Cycle Detection}

Given the forward assignment \( a^* : I \to S \) and the backward assignment \( b^* : S \to I \), we represent the resulting solution as a directed graph \( G = (V, E) \), where
\[
V = I \cup S,
\]
and edges alternate between pickups and deliveries:
\[
E = \{ (p_i \to s_{a^*(i)}) \mid i \in I \} \cup \{ (s_j \to p_{b^*(j)}) \mid j \in S \}.
\]

We apply a \emph{cycle detection algorithm} (e.g., depth-first search) to find all strongly connected components forming disjoint cycles:
\[
\mathcal{C} = \{ C_1, C_2, \ldots, C_m \}, \quad m \geq 1.
\]

If \( m = 1 \), the solution is a single cycle and the procedure terminates.

\paragraph{Step 2: Cycle Merging}

If multiple cycles exist, we iteratively merge pairs of cycles by modifying the edge structure to create a larger cycle, reducing the number of cycles by one each iteration.

\subparagraph{Cost Increment Evaluation}

For each pair of distinct cycles \( (C_a, C_b) \in \mathcal{C} \times \mathcal{C} \), we evaluate the \emph{incremental travel cost} \(\Delta T(C_a, C_b)\) incurred by merging them. This cost represents the increase in total travel distance when connecting \(C_a\) and \(C_b\) into a single cycle via a suitable edge swap.

Formally, let
\[
\Delta T(C_a, C_b) := \min_{\substack{(u,v) \in C_a \times C_b}} \big[ d(\text{end}(u), \text{start}(v)) - d(\text{old edges removed}) \big],
\]
where \(\text{end}(u)\) and \(\text{start}(v)\) denote nodes at which cycles \(C_a\) and \(C_b\) can be connected, and "old edges removed" correspond to edges broken to join the cycles. The exact calculation depends on the graph structure and allowable edge replacements. Details can be found in Appendix.

\subparagraph{Iterative Merging}

We select the pair \((C_{a^*}, C_{b^*})\) with the minimum \(\Delta T\) and perform the merge by appropriately redirecting edges to connect the two cycles into one. Then, update the cycle set:
\[
\mathcal{C} \leftarrow \mathcal{C} \setminus \{ C_{a^*}, C_{b^*} \} \cup \{ C_{merged} \}.
\]

This process repeats until only one cycle remains, ensuring the feasibility of the final tour.

\paragraph{Algorithmic Summary}

\begin{enumerate}
    \item Perform forward and backward assignments using the Hungarian algorithm.
    \item Construct the directed graph and detect all cycles.
    \item While the number of cycles \( m > 1 \):
    \begin{itemize}
        \item Evaluate \(\Delta T\) for all cycle pairs,
        \item Merge the pair with minimum \(\Delta T\),
        \item Update cycle set \(\mathcal{C}\) and cycle count \(m\).
    \end{itemize}
    \item Return the single merged cycle as the final routing solution.
\end{enumerate}

\begin{figure}[h!]
  \centering  \includegraphics[width=0.6\linewidth]{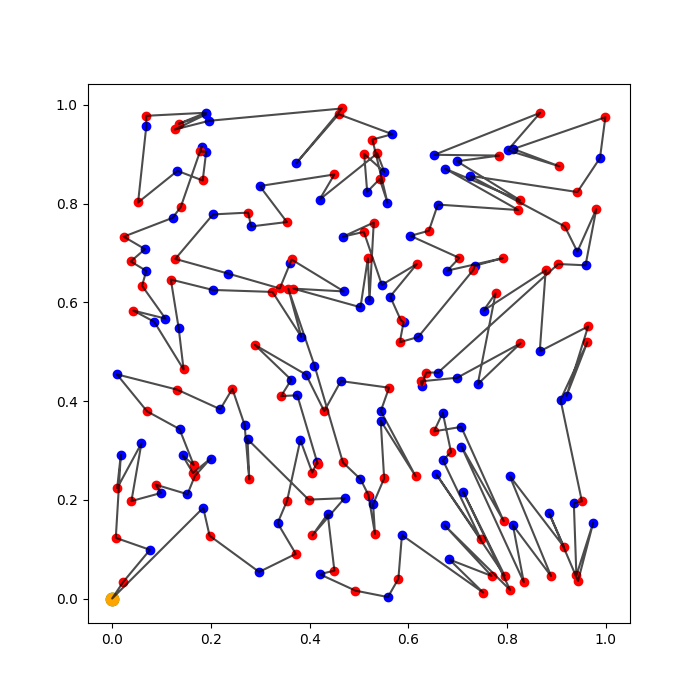}  
  \caption{Initial Guess from Cycle Merging Algorithm}
  \label{fig:f3}
  \vspace{0mm}
\end{figure}

This heuristic effectively balances minimizing additional travel cost with feasibility, yielding a complete single-cycle tour suitable for practical routing applications, which serves as a initial guess solution to be further improved, as shown in Figure \ref{fig:f3}.

\begin{figure}[h!]
  \centering  \includegraphics[width=0.6\linewidth]{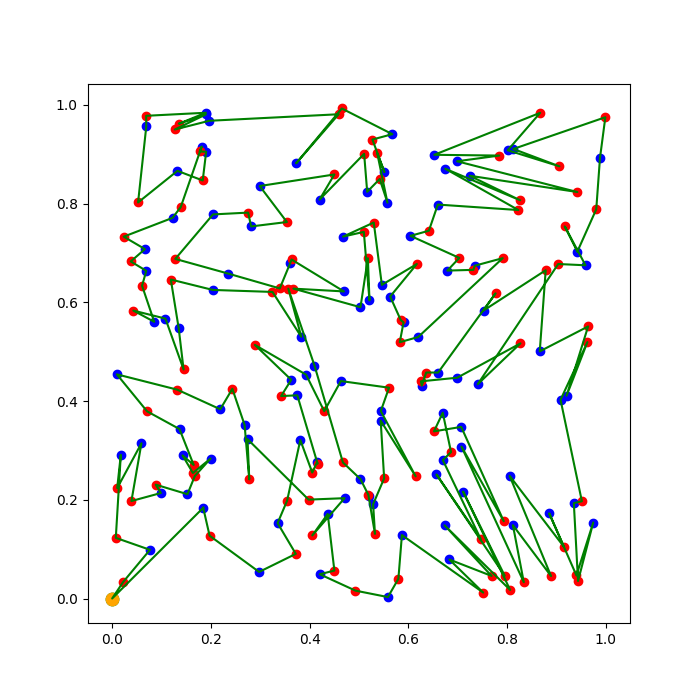}  
  \caption{Solution After Shaking and Simulated Annealing Process}
  \label{fig:f4}
  \vspace{0mm}
\end{figure}

\subsection{Shaking Method with Two-Way Cost Assignment}

To further improving the solution, we propose an iterative shaking method consisting of two alternating optimization steps that update the pickup sequence \(\pi\) and the assignment mapping \(a\) alternatively.

\subsubsection{Update Assignment \(a\) for Fixed Pickup Sequence \(\pi\)}

Given a fixed pickup sequence \(\pi = (\pi(0), \pi(1), \ldots, \pi(n))\), the placeholders are positioned between consecutive items in the sequence. For each position \(k = 0, \ldots, n\), we consider assigning a placeholder \(S_j\) between item \(I_{\pi(k)}\) and \(I_{\pi(k+1 \bmod (n+1))}\).

Define the cost matrix \(C^{(1)} \in \mathbb{R}^{(n+1) \times (n+1)}\), where each entry \(C^{(1)}_{k,j}\) represents the cost of assigning placeholder \(S_j\) to serve the transition from \(I_{\pi(k)}\) to \(S_j\), and from \(S_j\) to \(I_{\pi(k+1 \bmod (n+1))}\):

\[
C^{(1)}_{k,j} = d\big(p_{\pi(k)}, s_j\big) + d\big(s_j, p_{\pi(k+1 \bmod (n+1))}\big).
\]

This cost reflects the two-way travel through placeholder \(S_j\) placed between two consecutive items. The Hungarian algorithm is applied to minimize the total cost of this pairing, resulting in an updated assignment \(a\) such that:

\[
a_n = \arg\min_{a} \sum_{k=0}^{n} C^{(1)}_{k, j=a(\pi(k))}
\]
with the constraint \(a(0) = 0\) and bijectivity preserved.

\subsubsection{Update Pickup Sequence \(\pi\) for Fixed Assignment \(a\)}

Conversely, fix the assignment \(a\) and consider the induced placeholder sequence \(\pi_s\). Let \(\pi_s\) denotes the induced placeholder sequence from assignment \(a\), such that \(S_{\pi_s(k)} = a(\pi(k))\).

Construct a cost matrix \(C^{(2)} = [C^{(2)}_{i,j}] \in \mathbb{R}^{(n+1) \times (n+1)}\) where, for each candidate position \(j\), compute \(C^{(2)}_{i,j}\) as the incremental travel cost incurred by visiting item \(I_i\) between placeholders \(S_{\pi_s(j-1)}\) and \(S_{\pi_s(j)}\).

Formally,
\[
C^{(2)}_{i,j} = d\big(s_{\pi_s(j-1) mod(n+1)}, p_i\big) + d\big(p_i, s_{\pi_s(j)}\big),
\]
with indices modulo \(n+1\).

Apply the Hungarian algorithm on \(C^{(2)}\) to find the pickup sequence permutation \(\pi^{new}\) minimizing the total cost:
\[
\pi_{n} = \arg \min_{\pi \text{ permutation}} \sum_{j=0}^{n} C^{(2)}_{\pi(j), j}
\]

\subsubsection{Shaking Algorithm}
The shaking algorithm iteratively run the two updates mentioned above to improve the solution till converge to a local minima or achieve allowable steps, as show in Algorithm \ref{alg:shaking}. 

\begin{algorithm}[H]
\caption{Iterative Shaking Procedure}
\label{alg:shaking}
\KwIn{Initial pickup sequence $\pi^{(0)}$ with $\pi^{(0)}(0) = 0$, max iterations $T_{\max}$, tolerance $\epsilon$}
\KwOut{Optimized pickup sequence $\pi$ and assignment $a$}

$t \gets 0$\;
Compute initial assignment $a^{(0)}$ using Hungarian on cost matrix $C^{(1)}$\;
$T^{(0)} \gets T(\pi^{(0)}, a^{(0)})$\;

\Repeat{
  $|T^{(t)} - T^{(t-1)}| < \epsilon$ \textbf{or} $t \geq T_{\max}$
}{
    $t \gets t+1$\;
    
    Construct cost matrix $C^{(1)}$ from $\pi^{(t-1)}$\;
    $a^{(t)} \gets \text{Hungarian}(C^{(1)})$\;

    Construct cost matrix $C^{(2)}$ from $a^{(t)}$\;
    $\pi^{(t)} \gets \text{Hungarian}(C^{(2)})$, enforce $\pi^{(t)}(0) = 0$\;

    $T^{(t)} \gets T(\pi^{(t)}, a^{(t)})$\;
}

\Return $(\pi^{(t)}, a^{(t)})$\;
\end{algorithm}

\section{Computational Complexity Analysis and Efficiency Improved Methods}
\subsection{Brute-Force Computational Complexity}

For \textbf{assignment permutations:}  
the number of possible one-to-one assignments (bijections) between the \(n\) flexible items and \(n\) placeholders is \(n!\). For \textbf{pickup order permutations:}  
considering the fixed item \(I_0\) at the start of the sequence, the number of permutations of the remaining \(n\) items for pickup order is also \(n!\).
For each candidate solution, the total travel cost involves summations over pickup and delivery sequences of length proportional to \(n\).

Therefore, the overall computational complexity of the brute-force method is:
\[
  \mathcal{O}\bigl(n \times (n!)^2\bigr) ,
\]
which is super-exponential and becomes computationally infeasible even for moderate values of \(n\).

\subsubsection{Computational Complexity of Proposed Method}
In this subsection, we analyze the computational complexity of the proposed method.
The primary computational components of the method are:

\begin{itemize}
  \item \textbf{Assignment Optimization via Hungarian Algorithm:}  
  The initial assignment step involves solving an assignment problem on \(n\) flexible items and \(n\) placeholders (excluding the fixed pair \(I_0 \mapsto S_0\)), resulting in an \(n \times n\) cost matrix. The Hungarian algorithm operates with complexity \(\mathcal{O}(n^3)\) on this matrix.  
  The subsequent back-connection assignment consider all \(n+1\) items and placeholders, yielding an \((n+1) \times (n+1)\) matrix and requiring \(\mathcal{O}((n+1)^3)\) operations.  
  Thus, each iteration requires approximately:
  \[
    \mathcal{O}(n^3) + \mathcal{O}((n+1)^3) = \mathcal{O}((n+1)^3)
  \]
  operations for the two assignment steps combined.

  \item \textbf{Cycle Detection and Merging:}  
  Cycle detection in the assignment graph is performed using graph traversal algorithms such as depth-first search, with complexity \(\mathcal{O}(n)\) per iteration.  
  Cycle merging involves evaluating potential merges, which is bounded by \(\mathcal{O}(n^4)\) complexity. For large size problem when n>100, such merging time is the dominant time consumer in the solver. Therefore, efficient implementation of the merging code is very critical.

  \item \textbf{Iterative Shaking Process:}  
  Let \(T\) denote the total number of iterations until convergence or stopping criteria are met.  
  The overall complexity is approximately:
  \[
    \mathcal{O}\bigl(T \times (n^3+(n+1)^3 )\bigr),
  \]
  where the dominant term remains the cubic cost from the Hungarian algorithm executions.
\end{itemize}

The proposed iterative method with polynomial complexity per iteration represents a substantial computational improvement over brute-force enumeration, enabling practical solutions for larger problem instances.

\subsection{Partial Shaking Optimization with Simulated Annealing}

To further enhance scalability for large-scale instances, we introduce a global refinement heuristic based on \emph{partial shaking with simulated annealing (SA)}. This method reduces the computational burden of global optimization by operating on overlapping subproblems and incorporating stochastic exploration to escape local minima.

\paragraph{Motivation.}
The Hungarian algorithm and full assignment recomputation incur a cubic cost \(\mathcal{O}(n^3)\), which becomes expensive for large \(n\). Instead of optimizing over the full set of \(n\) items and placeholders, we propose optimizing over a smaller window of size \(m < n\), which significantly lowers the per-step cost while preserving global improvement.

\paragraph{Algorithm Overview.}
Let \(\pi = (\pi(0), \pi(1), \dots, \pi(n))\) denote the current item sequence, and \(a: I \to S\) be the current assignment. The method proceeds as follows:

\begin{enumerate}
    \item \textbf{Window Selection:}  
    Select a contiguous subsequence of items \(\pi^{\text{sub}} = \big(\pi(i), \pi(i+1), \dots, \pi(i+m-1)\big)\), where \(m < n\) is the window size and \(i = 0, l, 2l, \dots\) with step size \(l < m\).
    
    \item \textbf{Subproblem Formation:}  
    Extract the corresponding placeholders \(S^{\text{sub}} = \{a(\pi(k)) \mid k = i, \dots, i+m-1\}\). This pair \((\pi^{\text{sub}}, S^{\text{sub}})\) defines a smaller subproblem of size \(m\).

    \item \textbf{Simulated Annealing with Shaking:}  
    Perform simulated annealing (SA) within the subproblem to explore permutations of  only items and perform shaking algorithm to improve the candidates. Each SA step includes:
    \begin{itemize}
        \item Randomly update the item sequence \(\pi^{\text{sub}}\) in neighborhood.
        \item Optimize placeholders assignment and item sequence in the subsets using the local shaking algorithm that minimizes sub-cost.
        \item Evaluate the local cost function, incorporating entry and exit links from/to the global tour.
        \item Accept or reject the new state using the Metropolis criterion, based on temperature \(T(t)\).
    \end{itemize}

    \item \textbf{Update Global Tour:}  
    After convergence of the SA within the subproblem, update the global sequence \(\pi^{sub}\) and assignment \(a^{sub}\) with the improved solution $\pi_n^{sub}, a_n^{sub}$.

    \item \textbf{Sliding Window:}  
    Slide the window forward by \(l\) and repeat Steps 1–4 until the entire sequence has been processed.

    \item \textbf{Multiple Passes (Optional):}  
    Optionally repeat the full process multiple times with a decreasing temperature schedule to refine the global solution.
\end{enumerate}

\paragraph{Computational Efficiency.}
\begin{itemize}
    \item Each subproblem involves at most \(\mathcal{O}(m^3)\) operations (from the Hungarian method).
    \item The number of windows is approximately \(\lceil n / l \rceil\).
    \item Total complexity per pass: \(\mathcal{O}\left( \frac{n}{l} \cdot m^3 \right)\), which is significantly lower than \(\mathcal{O}(n^3)\) when \(m \ll n\).
\end{itemize}

Take note that only sampling the permutation of items would seldom make progress to the solution, because the new sequence is very likely to break the initial high quality connection to assigned placeholders. It is the shaking algorithm that reassign the placeholders to local minimal again and enable the possible progress for finding better solutions. 
When $m=n$, the Simulated Annealing algorithm is running for global solution searching, which is feasible for small size problem. However, for large size problem, it is very difficult to make effective progress because the entire size of the entire permutation is too large.  
This hybrid approach provides an efficient global refinement mechanism, integrating local structure-aware updates with global optimization capacity. Detailed algorithm for the Shaking enhanced SA methods in as showin in ALgorithm \ref{alg:partial_shaking_sa}

\begin{algorithm}[h!]
\DontPrintSemicolon
\SetAlgoLined
\SetKwInOut{Input}{Input}
\SetKwInOut{Output}{Output}

\caption{Partial Shaking with Simulated Annealing (with Best State Tracking)}
\label{alg:partial_shaking_sa}

\Input{Current pickup sequence $\pi$, assignment $a$, window size $m$, step size $l$, initial temperature $T_0$, temperature $T_{\min}$, cooling rate $\alpha$, sample number $n_s$}
\Output{Refined pickup sequence $\pi$, assignment $a$}

\For{$i \gets 0$ \KwTo $n - m + l$ \KwBy $l$}{
    Extract sub-sequence: 
    $\pi^{\text{sub}} \gets (\pi(i), \dots, \pi(\min(i + m - 1, n)))$ \;
    Extract sub-assignment: 
    $a^{\text{sub}} \gets \{ a(\pi(k)) \mid k = i, \dots, \min(i + m - 1, n) \}$ \;

    Initialize temperature: $T \gets T_0$ \;
    Set current state: $(\pi^{\text{curr}}, a^{\text{curr}}) \gets (\pi^{\text{sub}}, a^{\text{sub}})$ \;
    Set best state: $(\pi^{\text{best}}, a^{\text{best}}) \gets (\pi^{\text{sub}}, a^{\text{sub}})$ \;
    Set best cost: $C^{\text{best}} \gets \texttt{Cost}(\pi^{\text{best}}, a^{\text{best}})$ \;

    \While{$T > T_{\min}$}{
        Update current: $(\pi^{\text{curr}}, a^{\text{curr}}) \gets (\pi^{\text{best}}, a^{\text{best}})$ \;

        \For{$n \gets 0$ \KwTo $n_s$}{
            $C^{\text{curr}} \gets \texttt{Cost}(\pi^{\text{curr}}, a^{\text{curr}})$ \;
            
            $\pi^{\text{curr}} \gets sample_{neighbour}(\pi^{\text{curr}})$ \;
            
            $(\pi^{\text{new}}, a^{\text{new}}) \gets \texttt{ShakingLocalOptimize}(\pi^{\text{curr}}, a^{\text{curr}})$ \;
            $C^{\text{new}} \gets \texttt{Cost}(\pi^{\text{new}}, a^{\text{new}})$ \;

            $\Delta C \gets C^{\text{new}} - C^{\text{curr}}$ \;

            \eIf{$\Delta C < 0$ \textbf{or} $\texttt{Random}(0,1) < \exp(-\Delta C / T)$}{
                $(\pi^{\text{curr}}, a^{\text{curr}}) \gets (\pi^{\text{new}}, a^{\text{new}})$ \;
                \If{$C^{\text{new}} < C^{\text{best}}$}{
                    Update best: $(\pi^{\text{best}}, a^{\text{best}}) \gets (\pi^{\text{new}}, a^{\text{new}})$ \;
                    $C^{\text{best}} \gets C^{\text{new}}$ \;
                }
            }{
                Reject and retain current state \;
            }
            Update temperature: $T \gets \alpha \cdot T$ \;
        }
    }

    Replace subsegment in $\pi$ and $a$ with best state $(\pi^{\text{best}}, a^{\text{best}})$ \;
}
\end{algorithm}

\section{Evaluation of performance of Or-Tool base Solver}

To evaluate the time efficiency and solution quality of the proposed shaking-based optimization method, we benchmark the method with ground truth provided by Gurobi solver as shared in the GitHub repository, which will be discussed in next Section. In this section, we compare the introduced solution through establishing a baseline solution using a two-stage heuristic framework using Google Or-Tool. This baseline serves as a point of comparison for runtime and as a reference for approximation quality.

\subsection{Establishment of a Or-Tool based Solver}

The benchmark approach consists of the following steps:

\begin{enumerate}

    \item \textbf{Initial Assignment via Hungarian Algorithm}:  as shown in Section 3.1.

    \item \textbf{Sequence Optimization using TSP Solver}: \\
    Once the assignment is fixed, the sequence in which the agent picks up the items is optimized by formulating a Traveling Salesman Problem (TSP) with picking and delivering tasks over the item positions and placeholder positions. The tour starts from the fixed item \( I_0 \) and includes all other items exactly once. The Google OR-Tools TSP solver is employed for this task, configured with the ``Guided Local Search'' metaheuristic, to efficiently compute a near-optimal visiting order.

    \item \textbf{Local Refinement via Shaking Algorithm (Benchmark Variant)}: \\
    To assess the improvement potential of localized search, the above baseline solution is further refined using the proposed shaking-based optimization as a post-processing step. The shaking algorithm is applied directly to the complete set of TSP-derived item sequence and its corresponding assignment. This benchmark variant allows comparison between the standalone baseline and the shaking-enhanced baseline, in terms of total cost reduction and runtime.
\end{enumerate}
This solver can provide fair results, but the quality is limited even when we increase the solver timeout.  Table \ref{tab:alt_results} presents 10 samples of the results with timeout=10s for Or-Tool, comparing with the exact optimal solution from Gurobi as shared in Datasets \cite{yuan2025}. It is also noticed that longer runtime does not improve the solution results significantly. The direct solution from Or-Tool tool with feed-in assignment from Hungarian algorithm is provided in column two. Column three is the results after two more round of Or-Tool optimization, each followed by the introduced shaking algorithm (in Section 3.5) improvement, totally costing 30 over seconds with average of 1.66\% error from ground truth. 

\begin{table}[h!]
\centering
\caption{Comparison of alternative methods (OR-Tools and shaking) with Gurobi optimal solutions}
\begin{tabular}{|c|c|c|c|}
\hline
\textbf{Exp No.} & \textbf{Or-Tool (\%)} & \textbf{Final (\%)} & \textbf{Gurobi Opt(m)} \\
\hline
1 & 20.298 (3.90\%) & 20.226 (3.54\%) & 19.536 \\
2 & 31.132 (1.29\%) & 30.979 (0.80\%) & 30.734 \\
3 & 26.682 (2.19\%) & 26.307 (0.75\%) & 26.111 \\
4 & 19.607 (2.68\%) & 19.532 (2.28\%) & 19.095 \\
5 & 19.028 (2.60\%) & 18.747 (1.08\%) & 18.547 \\
6 & 19.653 (1.98\%) & 19.529 (1.33\%) & 19.273 \\
7 & 21.406 (3.02\%) & 21.017 (1.15\%) & 20.778 \\
8 & 25.677 (2.16\%) & 25.490 (1.42\%) & 25.134 \\
9 & 20.395 (2.36\%) & 20.355 (2.15\%) & 19.926 \\
10 & 21.177 (2.51\%) & 21.087 (2.07\%) & 20.660 \\
\hline
\textbf{Avg. Error (\%)} & \textbf{2.47\%} & \textbf{1.66\%} & – \\
\hline
\end{tabular}
\label{tab:alt_results}
\end{table}

\section{Experimental and Results}
\subsection{Experimental Design}

The performance of the proposed hybrid metaheuristic—comprising a merging heuristic followed by a introduced shaking algorithm guided simulated annealing (SA) refinement step—was evaluated across benchmark instances of increasing size: $n=32$, $n=100$, $n=200$, and $n=300$. Each problem size was tested on 10 independent instances using both the proposed approach and the exact solver Gurobi. Gurobi’s solutions serve as the benchmark ground truth cost, together with the computational time.

For each instance, we recorded:
\begin{itemize}
    \item \textbf{Gurobi's optimal solution cost and runtime}.
    \item \textbf{Average and standard deviation} of solution cost from 10 runs of our heuristic.
    \item \textbf{Minimum solution cost} achieved across the runs.
    \item \textbf{Initial merging solution's path distance cost ($L_M$)}.
    \item \textbf{Execution times}: separated into merging and SA phases (the shaking phase time is negligible comparably), and compared to Gurobi's runtime.
\end{itemize}

All experiments were executed on the same computational environment, for n=100 and above, parallel programming is applied to speed up the process. The heuristic’s performance is assessed based on relative deviation from optimality and scalability in runtime.

\subsection{Results and Analysis}

Tables~\ref{tab:n32XQ_results_updated} through~\ref{tab:experiment_results_n300_updated} present the experimental findings. Each Table presents the results of 10 experimental runs comparing the performance of the proposed heuristic method against the exact solutions obtained using Gurobi in previously shared Datasets. Each row corresponds to a distinct problem instance.  The second column lists the optimal cost computed by Gurobi, serving as the baseline for comparison. The third column reports the average solution cost obtained by executing the proposed method 10 times, accompanied by the population standard deviation (STDEV.P). The relative percentage error of the average cost with respect to Gurobi’s optimal is provided in parentheses.  The fourth column shows the best (minimum) solution cost achieved among the 10 heuristic runs, along with its corresponding percentage deviation from the optimal. The fifth column (denoted $L_M$) indicates the path cost of the initial merging solution.

Execution times are reported in the final columns: the combined runtime for the merge and SA phases ($dt_M + dt_{SA}$), and the Gurobi runtime for solving the same instance to optimality. 

\begin{table}[h!]
\footnotesize
\centering
\caption{Comparison Results for $n=32$ XQ}
\begin{tabular}{|c|c|c|c|c|c|c|}
\hline
\textbf{Exp No.} & \textbf{Gurobi} & \textbf{Average ± Std (\%)} & \textbf{Min (\%)} & \textbf{$L_M$ (\%)} & \textbf{$dt_M+dt_{SA}$ (s)} & \textbf{$dt_{Gurobi}$ (s)} \\
\hline
1 & 8.981 & 9.050 ± 0.012 (0.77\%) & 9.017 (0.40\%) & 9.107 (1.40\%) & 0.05 + 1.76 & 0.1749 \\
2 & 11.131 & 11.165 ± 0.013 (0.31\%) & 11.160 (0.26\%) & 11.246 (1.03\%) & 0.04 + 1.81 & 0.0184 \\
3 & 12.543 & 12.698 ± 0.001 (1.23\%) & 12.697 (1.23\%) & 12.705 (1.29\%) & 0.05 + 1.87 & 0.2018 \\
4 & 10.624 & 10.842 ± 0.005 (2.06\%) & 10.837 (2.01\%) & 10.891 (2.51\%) & 0.05 + 1.77 & 0.0195 \\
5 & 10.940 & 10.940 ± 0.000 (0.00\%) & 10.940 (0.00\%) & 10.940 (0.00\%) & 0.05 + 1.79 & 0.0450 \\
6 & 8.600 & 8.727 ± 0.000 (1.48\%) & 8.727 (1.48\%) & 8.744 (1.68\%) & 0.05 + 1.83 & 0.0611 \\
7 & 10.822 & 11.017 ± 0.024 (1.80\%) & 10.987 (1.52\%) & 11.124 (2.80\%) & 0.05 + 1.82 & 0.0149 \\
8 & 9.305 & 9.452 ± 0.044 (1.58\%) & 9.371 (0.71\%) & 9.711 (4.37\%) & 0.06 + 1.85 & 0.1165 \\
9 & 10.695 & 10.697 ± 0.003 (0.02\%) & 10.695 (0.00\%) & 10.708 (0.12\%) & 0.05 + 1.84 & 0.0742 \\
10 & 11.408 & 11.438 ± 0.008 (0.26\%) & 11.429 (0.18\%) & 11.453 (0.39\%) & 0.06 + 1.84 & 0.0331 \\
\hline
\textbf{Avg. \%} & – & \textbf{0.85\%} & \textbf{0.74\%} & \textbf{1.20\%} & – & – \\
\textbf{Avg. Time} & – & – & – & – & \textbf{0.05 + 1.80} & \textbf{0.0758} \\
\hline
\end{tabular}
\label{tab:n32XQ_results_updated}
\end{table}

\begin{table}[h!]
\footnotesize
\centering
\caption{Results of experiments for $n=100$, comparing heuristic and Gurobi solutions with execution times}
\begin{tabular}{|c|c|c|c|c|c|c|}
\hline
\textbf{Exp No.} & \textbf{Gurobi} & \textbf{Average ± Std (\%)} & \textbf{Min (\%)} & \textbf{$L_M$ (\%)} & \textbf{$dt_M + dt_{SA}$ (s)} & \textbf{$dt_{Gurobi}$ (s)} \\
\hline
0 & 19.536 & 19.670 ± 0.003 (0.69\%) & 19.663 (0.65\%) & 19.757 (1.12\%) & 0.14 + 1.07 & 0.5980 \\
1 & 30.734 & 30.837 ± 0.004 (0.34\%) & 30.828 (0.31\%) & 30.899 (0.54\%) & 0.11 + 1.10 & 0.7149 \\
2 & 26.111 & 26.275 ± 0.015 (0.62\%) & 26.247 (0.53\%) & 26.523 (1.57\%) & 0.18 + 1.09 & 0.8715 \\
3 & 19.095 & 19.286 ± 0.000 (1.01\%) & 19.286 (1.01\%) & 19.298 (1.06\%) & 0.14 + 1.10 & 0.9859 \\
4 & 18.547 & 18.650 ± 0.004 (0.56\%) & 18.637 (0.49\%) & 18.712 (0.88\%) & 0.13 + 1.07 & 2.0169 \\
5 & 19.273 & 19.374 ± 0.000 (0.53\%) & 19.374 (0.53\%) & 19.403 (0.68\%) & 0.16 + 1.10 & 1.5554 \\
6 & 20.778 & 21.179 ± 0.031 (1.94\%) & 21.151 (1.79\%) & 21.394 (2.95\%) & 0.17 + 1.08 & 0.8501 \\
7 & 25.134 & 25.185 ± 0.001 (0.20\%) & 25.183 (0.19\%) & 25.203 (0.28\%) & 0.13 + 1.14 & 0.2294 \\
8 & 19.926 & 20.023 ± 0.000 (0.48\%) & 20.023 (0.48\%) & 20.052 (0.63\%) & 0.13 + 1.15 & 0.2498 \\
9 & 20.660 & 20.758 ± 0.000 (0.47\%) & 20.758 (0.47\%) & 20.924 (1.27\%) & 0.14 + 1.12 & 2.6301 \\
\hline
\textbf{Avg. Error\%} & – & \textbf{0.68\%} & \textbf{0.64\%} & \textbf{1.12\%} & – & – \\
\textbf{Avg. Time} & – & – & – & – & \textbf{0.14 + 1.10} & \textbf{0.9709} \\
\hline
\end{tabular}
\label{tab:experiment_results_n100_updated}
\end{table}

\begin{table}[h!]
\footnotesize
\centering
\caption{Results of experiments for $n=200$, comparing heuristic and Gurobi solutions with execution times}
\begin{tabular}{|c|c|c|c|c|c|c|}
\hline
\textbf{Exp No.} & \textbf{Gurobi} & \textbf{Average ± Std (\%)} & \textbf{Min (\%)} & \textbf{$L_M$ (\%)} & \textbf{$dt_M + dt_{SA}$ (s)} & \textbf{$dt_{Gurobi}$ (s)} \\
\hline
0 & 32.526 & 32.850 ± 0.015 (0.99\%) & 32.835 (0.94\%) & 32.976 (1.38\%) & 0.73 + 1.60 & 1.8281 \\
1 & 29.779 & 29.948 ± 0.008 (0.57\%) & 29.929 (0.50\%) & 30.064 (0.95\%) & 0.68 + 1.61 & 8.3149 \\
2 & 26.949 & 27.427 ± 0.005 (1.77\%) & 27.418 (1.62\%) & 27.517 (2.08\%) & 0.66 + 1.66 & 4.2116 \\
3 & 27.459 & 27.768 ± 0.015 (1.11\%) & 27.738 (1.01\%) & 27.925 (1.74\%) & 0.63 + 1.55 & 6.2985 \\
4 & 28.546 & 28.901 ± 0.016 (1.23\%) & 28.857 (1.09\%) & 28.975 (1.52\%) & 0.74 + 1.56 & 14.5061 \\
5 & 32.552 & 32.834 ± 0.006 (0.87\%) & 32.825 (0.85\%) & 32.915 (1.13\%) & 0.68 + 1.56 & 6.9940 \\
6 & 38.105 & 38.335 ± 0.009 (0.61\%) & 38.320 (0.57\%) & 38.368 (0.69\%) & 0.60 + 1.57 & 8.7085 \\
7 & 38.060 & 38.273 ± 0.007 (0.56\%) & 38.263 (0.53\%) & 38.418 (0.95\%) & 0.67 + 1.59 & 4.1589 \\
8 & 31.544 & 31.678 ± 0.015 (0.42\%) & 31.665 (0.39\%) & 31.852 (0.98\%) & 0.61 + 1.59 & 10.0047 \\
9 & 28.029 & 28.160 ± 0.006 (0.47\%) & 28.150 (0.43\%) & 28.231 (0.73\%) & 0.57 + 1.55 & 11.0655 \\
\hline
\textbf{Avg. Error \%} & – & \textbf{0.87\%} & \textbf{0.77\%} & \textbf{1.22\%} & – & – \\
\textbf{Avg. Time} & – & – & – & – & \textbf{0.65 + 1.58} & \textbf{7.6993} \\
\hline
\end{tabular}
\label{tab:experiment_results_n200_updated}
\end{table}

\begin{table}[h!]
\footnotesize
\centering
\caption{Results of experiments for $n=300$, comparing heuristic and Gurobi solutions with execution times}
\begin{tabular}{|c|c|c|c|c|c|c|}
\hline
\textbf{Exp No.} & \textbf{Gurobi} & \textbf{Average ± Std (\%)} & \textbf{Min (\%)} & \textbf{$L_M$ (\%)} & \textbf{$dt_M + dt_{SA}$ (s)} & \textbf{$dt_{Gurobi}$ (s)} \\
\hline
0 & 42.002 & 42.324 ± 0.012 (0.77\%) & 42.314 (0.73\%) & 42.458 (1.07\%) & 1.75 + 2.08 & 76.0503 \\
1 & 33.129 & 33.634 ± 0.009 (1.54\%) & 33.623 (1.46\%) & 33.777 (1.94\%) & 2.20 + 2.09 & 25.6637 \\
2 & 34.331 & 34.652 ± 0.010 (0.93\%) & 34.642 (0.91\%) & 34.807 (1.37\%) & 1.78 + 2.10 & 44.0879 \\
3 & 35.198 & 35.503 ± 0.009 (0.87\%) & 35.488 (0.83\%) & 35.647 (1.28\%) & 1.72 + 2.09 & 28.5358 \\
4 & 33.812 & 34.037 ± 0.008 (0.67\%) & 34.024 (0.63\%) & 34.206 (1.14\%) & 1.62 + 2.10 & 77.8469 \\
5 & 37.478 & 37.706 ± 0.013 (0.61\%) & 37.696 (0.58\%) & 37.816 (0.91\%) & 1.66 + 2.11 & 84.7612 \\
6 & 33.954 & 34.396 ± 0.013 (1.29\%) & 34.391 (1.28\%) & 34.565 (1.76\%) & 1.72 + 2.09 & 26.5168 \\
7 & 40.668 & 40.806 ± 0.001 (0.34\%) & 40.804 (0.33\%) & 40.908 (0.59\%) & 1.78 + 2.11 & 28.3229 \\
8 & 36.677 & 36.956 ± 0.005 (0.76\%) & 36.947 (0.74\%) & 37.141 (1.27\%) & 1.91 + 2.11 & 51.3296 \\
9 & 38.658 & 38.859 ± 0.023 (0.52\%) & 38.820 (0.42\%) & 38.947 (0.74\%) & 1.76 + 2.10 & 83.7323 \\
\hline
\textbf{Avg. Error \%} & – & \textbf{0.85\%} & \textbf{0.74\%} & \textbf{1.20\%} & – & – \\
\textbf{Avg. Time} & – & – & – & – & \textbf{1.79 + 2.10} & \textbf{52.2849} \\
\hline
\end{tabular}
\label{tab:experiment_results_n300_updated}
\end{table}

The results clearly highlight the trade-off between optimality and runtime across problem sizes.

The proposed method consistently achieves near-optimal results. As the problem size increases, the accuracy remains remarkably stable. This stability across scales demonstrates consistency in generalization and robustness of the method. Moreover, the low standard deviations indicate solution consistency over multiple runs—critical for real-world reliability.

\subsection{Computational Efficiency}

The hybrid method significantly outperforms Gurobi in terms of runtime, especially as the problem size increases. While Gurobi’s exact optimization time scales superlinearly, the proposed heuristic maintains computational efficiency due to the lightweight merge and simulated annealing (SA) stages.

\begin{itemize}
    \item At $n=32$, Gurobi provides exact solutions in approximately \textbf{0.076 seconds}, while the heuristic remains similarly fast with negligible overhead.
    \item At $n=100$, the heuristic remains under \textbf{1.3 seconds}, whereas Gurobi averages around \textbf{0.97 seconds}. Although both are efficient at this scale, the heuristic’s flexibility allows consistent performance as the problem size grows.
    \item At $n=200$, the difference becomes evident: Gurobi requires about \textbf{7.7 seconds}, while the heuristic takes only \textbf{2.23 seconds} on average — a \textbf{3.5× speedup}.
    \item For $n=300$, Gurobi’s runtime increases sharply to \textbf{52.3 seconds}, while the heuristic completes in approximately \textbf{3.89 seconds}, achieving a \textbf{13× speedup}.
\end{itemize}

The heuristic design enables effective exploitation of problem structure (via merging), stochastic improvement (via shaking), and local refinement (via SA), all without incurring the exponential computational cost associated with combinatorial enumeration in exact solvers.

Additional experiments for $n=1000$ (available in the GitHub repository) show that the proposed hybrid method solves such large-scale instances in about \textbf{1 minute}, whereas Gurobi struggles to find feasible solutions within practical time limits.

Overall, the results confirm that the proposed approach not only delivers near-optimal solutions but also outperforms both Gurobi and OR-Tools in terms of \textbf{efficiency} and \textbf{scalability}, making it well-suited for solving larger and more complex optimization problems.

\subsubsection{Comparability with Other Algorithms}

A key advantage of this approach is its modularity. Each component—assignment, sequence optimization, cycle merging, local search—can be independently enhanced or replaced. This design allows seamless integration of additional constraints such as 2-Opt, 3-Opt, k-Opt, LKH algorithms etc. Furthermore, the shaking procedure can be extended into a broader metaheuristic framework, such as simulated annealing, variable neighborhood search, or large neighborhood search, enabling even more exploration of the solution space.

\subsubsection{Limitations and Future Work}
The cycle merging algorithm works better than a naive greedy algorithm, this algorithm only consider merging two cycle at a time. Therefore, there are still spaces to further improve the path merging quality using more advance methods. The shaking algorithm has its limitation that it is always updating the assignment and the sequence separately which only converges to local minima. For SA algorithms, it would be very time consuming to explore huge number of samples and thus it mainly works for local solution improvements. For global solution optimization, further research need to be conducted to improve the solution.

\section{Conclusion}
\label{sec:conclusion}

This work focuses on solving the joint assignment and routing optimization problem. While exact solvers such as Gurobi can provide optimal solutions for small to moderate problem sizes—as demonstrated in our prior work—they are not fast enough to scale to larger instances commonly encountered in real-world applications.

To address this gap, we introduce a fast metaheuristic framework that achieves near-optimal results with significantly lower computational cost. Empirical evaluations show that the proposed method consistently produces solutions within 1\% of optimality, while outperforming general-purpose solvers like OR-Tools in both accuracy and speed.

Although the approach proves effective and scalable, it represents one step toward efficient large-scale solvers for joint assignment and routing tasks. The current approach does not explore in depth on global optimization heuristic solution updates, which limits the performance of the solution accuracy to exact solution.  Future research can build on this direction by exploring more adaptive heuristics, hybrid algorithms, or apply them in task and motion planning in robotic domains.
\section{Appendix}
\subsection{Incremental Cost Function for Cycle Merging}
\begin{algorithm}[H]
\caption{Cycle Detection and Merging}
Compute forward assignment \(a^*\) and backward assignment \(b^*\) using Hungarian algorithm\;
Construct directed graph \(G\) from \(a^*, b^*\)\;
Detect cycles \(\mathcal{C} = \{C_1, \dots, C_m\}\) in \(G\)\;
\While{\(m > 1\)}{
    Initialize \(\Delta T_{\min} \gets +\infty\)\;
    \ForEach{distinct pairs \((C_a, C_b) \in \mathcal{C}\)}{
        \ForEach{edges \((v_i^a, v_{i+1}^a)\) in \(C_a\)}{
            \ForEach{edges \((v_j^b, v_{j+1}^b)\) in \(C_b\)}{
                Compute \(\Delta T_{i,j}\) as above\;
                \If{\(\Delta T_{i,j} < \Delta T_{\min}\)}{
                    Update \(\Delta T_{\min} \gets \Delta T_{i,j}\)\;
                    Record best merge edges \((v_i^a,v_{i+1}^a), (v_j^b,v_{j+1}^b)\)\;
                    Record cycles \((C_a, C_b)\)\;
                }
            }
        }
    }
    Perform merge by replacing recorded edges to connect \(C_a\) and \(C_b\)\;
    Update \(\mathcal{C}\) by removing \(C_a, C_b\) and adding merged cycle\;
    Update \(m \gets m - 1\)\;
}
Return single merged cycle \(C\)\;
\end{algorithm}
Merging cycles of sub-tours is very critical in resulting initial guess solution with great quality. In general, multiple cycles can be merged or together with more optimal choices, however, it would be more challenging and computationally more costly. 
Here, we introduce algorithm to merge two distinct cycles. 
To merge two cycles \(C_a\) and \(C_b\) into a single cycle, we remove one edge from each cycle and add two new edges to connect them, maintaining the alternating pickup-delivery structure and feasibility.

Let
\[
C_a = (v_1^a, v_2^a, \dots, v_{k_a}^a), \quad
C_b = (v_1^b, v_2^b, \dots, v_{k_b}^b)
\]
be the node sequences of cycles \(C_a\) and \(C_b\).

\paragraph{Edge Replacement}

We consider breaking edges
\[
(v_i^a \to v_{i+1}^a) \in C_a, \quad
(v_j^b \to v_{j+1}^b) \in C_b,
\]
and replacing them with cross edges while ensuring the item-placeholder connection (in such problem only one replacing connection is available). 
\[
(v_i^a \to v_{j+1}^b), \quad (v_j^b \to v_{i+1}^a).
\]

This operation merges \(C_a\) and \(C_b\) into a single cycle by reconnecting the node sequences.

\paragraph{Incremental Cost}

The incremental cost \(\Delta T_{i,j}\) for this merge is given by
\[
\Delta T_{i,j} = d(v_i^a, v_{j+1}^b) + d(v_j^b, v_{i+1}^a) - d(v_i^a, v_{i+1}^a) - d(v_j^b, v_{j+1}^b).
\]

We evaluate \(\Delta T_{i,j}\) for all valid edge pairs \((i,j)\) and select the pair with minimal \(\Delta T_{i,j}\).

The pair \((C_a, C_b)\) with minimal \(\min_{i,j} \Delta T_{i,j}\) across all cycle pairs is selected for merging.
\[
(i^*, j^*) = \arg\min_{i,j} \Delta T_{i,j}.
\]

\bibliographystyle{unsrtnat}
\bibliography{references} 
\end{document}